\documentclass[12pt,draft]{amsart}
\usepackage{amsmath,amsfonts,amsthm}
\usepackage[all]{xy}
\usepackage{fullpage}
\usepackage{amssymb}
\usepackage{colordvi}
\usepackage{url}

\numberwithin{equation}{section}

\newtheorem{theorem}{Theorem}[section]
\newtheorem{proposition}{Proposition}[section]

\newtheorem{corollary}{Corollary}[section]

\theoremstyle{definition}
\newtheorem{definition}{Definition}[section]

\theoremstyle{remark}
\newtheorem{remark}{Remark}[section]
\newtheorem{example}{Example}[section]

\newcommand{\CC}{\mathbb{C}}
\newcommand{\FF}{\mathbb{F}}

\newcommand{\QQ}{\mathbb{Q}}

\newcommand{\ZZ}{\mathbb{Z}}

\newcommand{\cA}{\mathcal{A}}
\newcommand{\cB}{\mathcal{B}}

\newcommand{\cO}{\mathcal{O}}
\newcommand{\cI}{\mathcal{I}}

\newcommand{\gC}{\mathfrak{C}}
\newcommand{\gP}{\mathfrak{P}}
\newcommand{\gR}{\mathfrak{R}}
\renewcommand{\aa}{\mathfrak{a}}

\newcommand{\pp}{\mathfrak{p}}
\newcommand{\qq}{\mathfrak{q}}

\newcommand{\lc}{\mbox{leadcoeff}}
\newcommand{\isom}{\cong}

\newcommand{\clos}{\overline} 
\newcommand{\conj}{\overline} 
\newcommand{\redu}{\overline} 

\renewcommand{\>}{\rangle}

\newcommand{\Aut}{\mathrm{Aut}}
\newcommand{\End}{\mathrm{End}}
\newcommand{\Gal}{\mathrm{Gal}}

\newcommand{\Jac}{\mathrm{Jac}}

\newcommand{\ignore}[1]{}

\begin{document}

\bibliographystyle{alpha}

\title[The $p$-adic CM-method for genus 2]{The $p$-adic CM-method for genus 2}
\author{P. Gaudry, T. Houtmann, D. Kohel, C. Ritzenthaler, and A. Weng}
\address{Laboratoire d'Informatique (LIX)\\
{\'E}cole polytechnique\\
91128 Palaiseau CEDEX\\
France
}
\email{\begin{tabular}{l}gaudry,houtmann,weng@lix.polytechnique.fr,kohel@maths.usyd.edu.au,\\ritzenth@math.jussieu.fr\end{tabular}}
\thanks{The second author thanks the ARC for financial support.\\\indent The fourth
  author thanks the mathematical departement of Sydney university for its warm
  hospitality.\\\indent The fifth author thanks the DFG for financial support.}
\date{\today}

\begin{abstract}
We present a nonarchimedian method to construct hyperelliptic
CM-curves of genus 2 over finite prime fields.
\end{abstract}
\maketitle

\begin{center} \begin{minipage}{15cm} \begin{small}
Throughout the document we use the following conventions (this is only
for the reference and use of the authors):
\begin{center}
\begin{tabular}{ll}
$d$ & degree of the base field of the curve, i.e.~$C/\FF_{2^d}$\\
$s$ & number of isomorphism classes, in elliptic curve case $s = h_K$\\
$n$ & degree of an irreducible component of class invariants\\
$K$ & a CM field\\
$K_0$ & the real subfield of $K$\\
$K^*$ & the reflex CM field of $K$\\
$K_0^*$ & the real subfield of $K^*$\\
$j_1$ & absolute Igusa invariant $J_2^5J_{10}^{-1}$\\
$j_2$ & absolute Igusa invariant $J_2^3J_4J_{10}^{-1}$\\
$j_3$ & absolute Igusa invariant $J_2^2J_6J_{10}^{-1}$\\
$N$ & $2$-adic precision\\
\end{tabular}
\end{center}
\end{small} \end{minipage} \end{center}

\section{Introduction}
\label{s:intro}

In 1991 Atkin proposed an algorithm for constructing elliptic curves over
finite fields with a given endomorphism ring~\cite{AtkinA,AtkinB}.
This algorithm originally proposed to speed up the Goldwasser-Kilian
primality test has several applications. Since the knowledge of the
endomorphism ring, enables us to easily determine the number of points
on the elliptic curve, it can for example be used to construct elliptic
curves with a prime order which has applications to cryptography.
The complex multiplication method has also become attractive to construct
suitable curves for pairing based cryptography~\cite{dem,barreto,BW}.

The usual CM-method works with floating point arithmetic.
We first construct all $h = h(\cO)$ isomorphism classes of elliptic
curves with complex multiplication by a given order $\cO$ of
discriminant $D$ in an imaginary quadratic field $K = \QQ(\sqrt{D})$.
We then compute their $j$-invariants numerically and build the minimal
polynomial
$$
H_D(X)=\prod_{i=1}^{h} (X-j_i)
$$
which by theory has integer coefficients, that can be recognized from
their floating point value if the precision of the computation is high
enough.
The CM-method has been generalized to higher genus, i.e.~to genus 2
curves and some special cases in genus~3~\cite{WengA,WengB,WengC}.

Recently, nonarchimedian approaches to the construction of class
polynomials $H_D(X)$ and analogues have been developed~(see~\cite{CH,BS}).
In this setting, given a imaginary quadratic order $\cO$ of discriminant
$D$ we choose a prime $p$ of size roughly $D$ such that there exists an
elliptic curve with complex multiplication by $\cO$ over $\FF_p$ ---
such a curve is found by exhaustive search.
A canonical lift of the $j$-invariant of the initial curve is computed
$p$-adically to sufficient precision to recover its minimal polynomial
$H_D(X)$ over $\QQ$.

In this paper we consider an analogue of the nonarchimedian approach
to construct class polynomials of hyperelliptic curves of genus 2.
We use a higher dimensional generalizations of the AGM over a $2$-adic
field.

Our paper is organized as follows. We first demonstrate the basic idea
by using an example in genus 1 (see Section \ref{genus1}).
We then recall some theoretical facts on complex multiplication of
abelian varieties of dimension 2 (see Section \ref{cmfacts}).
In Section \ref{main} we give an overview of the complete algorithm.

For our algorithm we need to explain how to run over isomorphism classes
of ordinary genus $2$ curves in characteristic $2$ (see Section
\ref{enum}). We need also to describe the AGM method for hyperelliptic
curves of genus 2 (see Section \ref{agm2}). We revise the $p$-adic
$LLL$-algorithm and describe some modifications which are specific to
our situation (see Sections \ref{lll}).
We also discuss how to determine the endomorphism ring of
a hyperelliptic curve over characteristic 2 in special cases
(Section \ref{maxorder}).

Finally, we give numerical examples which show that the $p$-adic
method can be efficiently used to compute class polynomials of
certain quartic CM fields (see Section \ref{ex}).

\section{Description of the basic AGM method for elliptic curves}
\label{genus1}

We first recall the AGM method for elliptic curves and explain
how it can be used to generate the class polynomial for imaginary
quadratic fields $K=\QQ(\sqrt{D})$ with $D\equiv 1\bmod 8$.
Let $k$ be a $2$-adic local field with uniformizer $\pi$ and let $a,b\in k$
be two elements such that
$$
\frac{a}{b}\equiv 1\bmod (8\pi).
$$
We can then take the square root $x=\sqrt{a/b\,}$ which is uniquely
determined if we impose the condition $x\equiv 1\bmod 4\pi$.
The sequence of pairs defined by
\begin{align*}
(a_{i+1},b_{i+1}) = (\frac{a_i+b_i}{2}, b_i\sqrt{\frac{a_i}{b_i}})
\end{align*}
derive from $2$-isogenies between elliptic curves.  More precisely,
if $E$ is an elliptic curve given by an equation of the form
$$
E_i: y^2 = x(x-a_i^2)(x-b_i^2),
$$
then the curve
$$
E_{i+1}: y^2=x(x-a^2_{i+1})(x-b^2_{i+1})
$$
is 2-isogenous to $E_i$ (possibly over some extension).
Moreover the value $t_i = a_i/b_i$ is an isomorphism invariant of
the pair $(E_i,E_{i+1})$ with their full $2$-torsion structures,
and if $E_i$ is defined over the unramified extension of $\QQ_2$,
then $E_i \rightarrow E_{i+1}$ reduces to the Frobenius modulo $2$.

Suppose that we are given an ordinary elliptic curve $E$ over $\FF_q$ with $q=2^d$
for some $d$. Let $\QQ_q$ be the unique unramified extension of degree $d$ of
$\QQ_2$. Then by a Theorem of J. Lubin, J.-P. Serre and J. Tate
\cite{lst}, there exists an elliptic curve $\tilde{E}/K_q$ such that
$$
\End(E)\simeq \End(\tilde{E}).
$$
This curve is called the canonical lift of $E$.
Given $E$ with $\End(E)\simeq \cO_K$ we want to construct the
polynomial $H_D(X)$.
Suppose that $h_K=t\cdot d$, i.e.~the prime $\pp_2$ lying above
$2$ in $\QQ(\sqrt{D})$ has order $d$ in the class group.
We use the AGM method to construct a cycle of 2-isogenous elliptic curve
$$
E_0^1 = E \rightarrow
E_1^1 \rightarrow
\cdots \rightarrow E_{d}^{1}=E_0^2
$$
where $\varphi_i:E_i\rightarrow E_{i+1}$ is a 2-isogeny.
Repeating the cycle sufficiently many times, we get a sequence of
elliptic curves such that $E_{j}^i$ is a good approximation for
the canonical lift of $E_j$.
We can then recover the $j$-invariant of $\tilde{E}_j$ with
high precision. If $d\ne h_K$,
we have to repeat this process for other elliptic curves of $\FF_q$ with the
same endomorphism ring until we found all $j$-invariants in $\QQ_q$. We
then compute the class polynomial with coefficients in $\QQ_q$ that we
recognize as integers if the precision is high enough.
\begin{example}
Consider a simple example. Let $D=-15$. Then $2$ splits into two
nonprincipal prime ideals in $K = \QQ(\sqrt{D})$ and we find the curve
$$
E: y^2+ xy = x^3 + \alpha^2
$$
over $\FF_2(\alpha)=\FF_2[x]/(x^2+x+1)$ with
$\End(E)=\cO_K=\ZZ[(1+\sqrt{-15})/2]$.
We lift $E$ to the curve $\tilde{E}/K_q$ where $\QQ_q$ is
the unramified extension of $\QQ_2$ of degree $2$ given by
$$
\tilde{E}: y^2=x(x-a_0^2)(x-b_0^2)
$$
where $\beta$ is a lift of $\alpha$ to $K_q$ and $a_0=1+4\beta^2$ and
$b_0=1-4\beta^2$. We now apply 13 rounds of the AGM and obtain
\begin{align*}
j(E_{1}^{13})=8026247402149799202321\beta - 6102896026815785332240,\\
j(E_{2}^{13})=3730718496258231955951\beta + 2950325125578927178719.
\end{align*}
The Hilbert class polynomial, determined modulo $2^{28}$, is given
by
$$
H_{-15}(X) = X^2 + 191025X - 121287375.
$$
\end{example}

\noindent{\bf N.B.}
The size the coefficients of $H_D(X)$ can be explicitly bounded by
$$
\frac{\pi\sqrt{|D|}}{\ln(10)}\sum \frac{1}{a}+10,
$$
where the sum runs over all reduced quadratic forms $(a,b,c)$ of
discriminant $D$ (see~\cite[p.~416]{cohen}), so the precision needed
for this algorithm can be effectively determined.

There are two main obstructions to extending to an arbitrary
discriminant~$D$.  First, the size of the coefficients of the
output polynomial $H_D(X)$ makes the construction of the Hilbert
class polynomial expensive even for $D$ of modest size, and
second, the application of the AGM imposed a congruence condition
$D \equiv 1 \bmod 8$.
In order to achieve a reduction in the coefficients size, one
can use alternative modular functions, e.g.~on some modular curve
$X_0(N)$.
In the AGM example, the modular invariant $t_i = a_i/b_i$ is a
function on $X_0(8)$ of the form $(u_i+4)/(u_i-4)$, and the AGM
recursion determines a lifted invariant which satisfies the
smaller minimal polynomial
$$
X^4 - 9X^3 + 17X^2 + 24X + 16,
$$
any root $u$ of which determines a CM $j$-invariant by
$$
j = \frac{(u^4 + 224u^2 + 256)^3}{u^2(u+4)^4(x-4)^4}.
$$
Existence of generalised AGM methods for elliptic curves in odd
characteristic have been proved by Carls \cite{Carls}.  Explicit
formulae for AGM recursions described for modular functions on
various $X_0(pN)$ and small characteristics $p$ were determined
by Kohel~\cite{AGMX0N} and by Br\"oker and Stevenhagen~\cite{BS}
using Weber functions (modular functions of level $N = 48$) and
small characteristic~$p$.

The method of Couveignes and Henocq~\cite{CH} for level $N = 1$
imposes no congruence condition on the input discriminant,
while variation of the level $N$ also varies of the congruence
condition on~$D$.  In another direction, Lercier and
Riboulet-Deyris~\cite{LerRibD04} use a $p$-adic lift of a CM
order embedded in the endomorphism ring of a supersingular
elliptic curve.  For $p = 2$ this allows one to treat the
complementary classes $D \equiv 5 \bmod 8$ and (fundamental)
$D \equiv 0 \bmod 4$.

\begin{remark}
\label{remell}
The case of genus 1 can be compared and contrasted with the
problems which arise in the generalisation to higher dimension.

1. To determine the class polynomial of a maximal order $\cO$,
we have to ensure that a selected curve $E/\FF_q$ has complex
multiplication by $\cO$ and not some suborder.
Determining the correct order $\End(E)$ requires a more detailed
analysis~(for example see~\cite{kohel} and some extensions to
genus $2$ in~\cite{EisLau04}).

2. For an elliptic curve $E/\FF_q$, such that its $j$-invariant
generates $\FF_q$ over $\FF_p$, the class number must be divisible
by the extension degree $d = [\FF_q:\FF_p]$.  The order $\End(E)$
of a randomly chosen $E$, however, has discriminant $D = t^2-4q$,
whose class number tends to grow like $O(\sqrt{q})$. In the case of genus
2, the class number will tend to grow faster.

3. All elliptic curves over a finite field $\FF_q$ which have
complex multiplication by an imaginary quadratic order $\cO$
have the same field of definition.
This follows from the Galois theory of class fields for
imaginary quadratic fields; its generalization to higher dimension
does not preserve this feature.

4. The $j$-invariant is an algebraic integer and we have explicit
bounds on the size of $j$ in terms of the discriminant of the order.
The lack of explicit bounds and the failure of the Igusa invariants
to be algebraic integers provide both technical and theoretical
obstacles.  As a result, even proving the correctness of the result
becomes more cumbersome (see Section~\ref{ex}).
\end{remark}

\section{The theoretic background}
\label{cmfacts}

In this section we will summarize some basic facts on Jacobians of
genus two curves and quartic CM fields needed to understand the
algorithm represented in the next section.

\subsection{The Frobenius endomorphism and its characteristic polynomial}

Let $C$ be a hyperelliptic curve of genus 2 over a finite field $\FF_q$ and
let $J_C$ be the Jacobian of $C$. Note that $J_C$ is an abelian surface.
Let $\pi_q$ be the Frobenius
endomorphism on $J_C$. Let $T_\ell$ be the Tate module for $J_C$ for some
prime $\ell$, $(\ell,q)=1$. The
Frobenius operates on the $4$-dimensional vector space $T_\ell\otimes \QQ_\ell$
and the characteristic polynomial $f_{\pi_q}(x)\in\ZZ[x]$ of this
representation is independent of the prime $\ell$. It classifies
the isogeny class of the Jacobian over $\FF_q$.
Any root $w$ of the Frobenius polynomial has absolute value
$\sqrt{q}$ and we have
\begin{equation}\label{grord}
f_{\pi_q}(1)=\#J_C(\FF_q).
\end{equation}
Given $f_{\pi_q}(x)$, we can determine
$\End(J_C)\otimes\QQ$~\cite{Waterhouse}.
If $f_{\pi_q}(x)$ is irreducible, then the Frobenius endomorphism
generates a CM field of degree 4, i.e.~a totally imaginary quadratic
extension of a real quadratic field.
\medskip

In the opposite direction, in Section \ref{main} we will try to
construct a curve over a large finite field $\FF_p$ whose Jacobian
has complex multiplication by the maximal order in a given CM
field $K$. Suppose we have given such
a curve $C$. The Frobenius endomorphism on $J_C$ corresponds to an
element $w\in\cO_K$ with absolute value $\sqrt{q}$. If we know
that $J_C$ is simple, then $\QQ(w)=K$. This will always be the case if
$K$ is non-normal or cyclic.

There are only finitely many elements in $\cO_K$
such that $w\conj{w}=q$. For each $w$ we can compute the minimal
polynomial
$f_w(x)\in\ZZ[x]$. If the Jacobian is ordinary and $K$ does
not contain any nontrivial roots of unity, we find precisely two
different $w$ up to conjugation and two different group orders in the
Galois case and two or four different $w$ up to conjugation and two,
three or four different group orders in the non-normal case
(cf.~\cite{WengH}).

We can now find the right order $n$ by choosing random
elements in the Jacobian and multiplying them with the possible values
for $n$.

\subsection{Quartic CM fields}\label{cmiso}
Let $K$ be a quartic CM field and let $\Phi=\{\varphi_1,\varphi_2\}$
be a set of two different embeddings of $K$ into $\CC$ such that
$\varphi_1\ne \varphi_2\rho$ where $\rho$ is the complex conjugation.
Then $(K,\Phi)$ is called a CM type; up to conjugation, there exist
exactly two different CM types.
To every abelian variety over $\CC$ with complex multiplication by
an order of $K$ we can assign a specific CM type.  This CM type
is called {\it primitive} if and only if the abelian variety is
absolutely simple.
A quartic CM field may be non-normal (whose normal closure is a $D_4$
extension of $\QQ$), cyclic, or bicyclic Galois extensions of $\QQ$.
For the first two, every CM type is primitive, but every bicyclic CM
type is nonprimitive, so we focus on the case that $K$ is non-normal
or cyclic over $\QQ$.

We can show that conjugate CM types will lead to the same set of
isomorphism classes of abelian varieties.
In the cyclic case, the set of isomorphism classes of one specific
CM type coincides with the set of isomorphism classes of any other
CM type.
Hence it will be enough to consider only one fixed CM type
(cf.~\cite{Spall1}).
For a CM field $K$, we denote by $\cO_K$ its maximal order and by $K_0$
its quadratic real subfield.
In order to
determine the number $s$ of isomorphism classes of principally
polarized abelian varieties with CM by $\cO_K$, we define an
associated class group.

\begin{definition}
Let $\cI(K)$ be the group of fractional ideals in $K$, and let
$K^\times$ act on the group $\cI(K) \times K_0^\times$ by
$\mu(\aa,\alpha) = (\mu\aa,\alpha\mu\conj{\mu})$.
Then the subgroup of $\cI(K) \times K_0^\times$ consisting
of pairs $(\aa,\alpha)$ such that $\aa\conj\aa = (\alpha)$
for totally positive $\alpha \in K_0$ contains the image of
$K^\times$, and we define $\gC(\cO_K)$ to be quotient of this
subgroup by $K^\times$.
\end{definition}

\noindent
The following theorem summarises the results of $\S14.6$ of
Shimura~\cite{ShimuraI}, and provides the explicit class number
for the set of isomorphism classes of principally polarised
CM abelian varieties.

\begin{theorem}
The set of isomorphism classes of principally polarised abelian
varieties with CM by $\cO_K$ is a principal homogeneous space
over $\gC(\cO_K)$, in particular $s = |\gC(\cO_K)|$.
\end{theorem}

We note that $\gC(\cO_K)$ is an extension by a group of order $1$
or $2$, of the kernel of the norm map $Cl(\cO_K) \rightarrow
Cl^+(\cO_{K_0})$
\footnote{
Here $Cl^+(\cO_{K_0})$ is the group of ideals modulo totally
positive principal ideals, and $Cl^+(\cO_{K_0}) = Cl(\cO_{K_0})$
if the fundamental unit of $\cO_{K_0}$ has norm $-1$.}
given by $\aa \mapsto \aa\conj\aa$.
The class of $(\cO_K,1)$ is the identity element, and $\gC(\cO_K)$
is a nontrivial extension of this kernel if and only if the
fundamental unit $\epsilon_0$ of $\cO_{K_0}$ has norm $1$ and
is not in the image of a fundamental unit of $\cO_K$.
In this case $(\cO_K,\epsilon_0)$ is a second element of
$\gC(\cO_K)$ which lies over the principal class of $Cl(\cO_K)$.

\begin{corollary}
Let $s$ denote the number of isomorphism classes of principally
polarised abelian varieties with CM by a maximal CM order $\cO_K$,
and let $h^\prime$ be the order of the kernel of the norm map
$Cl(\cO_K) \rightarrow Cl^+(\cO_{K_0})$.
If $N(\epsilon_0)$ equals $-1$, then $s = h^\prime$ if the field
$K$ is normal and $s = 2h^\prime$ if $K$ is non-normal.
If $N(\epsilon_0)$ equals $1$, then $s = h^\prime$ if $\epsilon_0$
is the norm of a unit in $\cO_K$, and $s = 2h^\prime$ otherwise.
\end{corollary}

The Cohen-Lenstra heuristics imply that the class number of the real
quadratic field $K_0$ has class number 1 with density greater that
$3/4$. In this case we can express a more precise form of the theorem
(see~\cite{WengH}).

\begin{corollary}
Let $K$ be quartic CM field, with real quadratic subfield $K_0$ of
class number~1.  If $K$ is cyclic over $\QQ$, then there are $h_K$
isomorphism classes, and if $K$ is not normal over $\QQ$ then there
are $2h_K$ isomorphism classes, with $h_K$ classes associated to
each CM type.
\end{corollary}

\noindent{\bf N.B.}
The enumeration of the isomorphism classes does not provide the
Galois action on their moduli.
The CM moduli determine an abelian extension of the Galois group
$Cl(\cO_{K^*})$ of the {\it reflex field} $K^*$ via a map
$Cl(\cO_{K^*}) \rightarrow \gC(\cO_K)$.
In the cyclic case, $K^*$ and $K$ coincide, but in the non-normal
case, $K$ and $K^*$ are nonisomorphic quartic CM fields embedded in
the normal closure $L$ of $K$.  In the latter case, the action on
the CM isomorphism classes is given by $\aa \mapsto (g(\aa),N(\aa))$
where $g$ is the composition of ideal extension to $L$ with the
norm of $L/K$, and $N = N_{K^*/\QQ}$.
Even if $Cl(\cO_{K^*}) \isom Cl(\cO_{K})$, the map to $\gC(\cO_K)$
may have a kernel which results in reducibility of the corresponding
class equations (see
Shimura~\cite[Main Theorem~1, Note~3, pp.~112-113]{ShimuraI}).

\subsection{The splitting of a prime in a given CM field}
\label{split}
Let $K$ be a quartic CM field. Analogously to the elliptic curve case
we can define invariants which classify the isomorphism class of the
hyperelliptic curves of genus $2$ or equivalently the principally
polarized abelian surfaces over $\CC$ completely. In contrast to the elliptic curve
case, the moduli space is $3$-dimensional and we find three
$j$-invariants $j_1$, $j_2$, $j_3$. We define the class polynomial
$$
H_k(X)=\prod_{\sigma\in\Sigma} (X-j_k^\sigma),\quad k=1,2,3
$$
where $\Sigma$ is the set of all isomorphism classes of principally
polarized abelian surfaces with $CM$ by the maximal order
$\cO_K$. Since we run over all isomorphism classes, the
polynomials $H_k(X)$ are Galois invariants, i.e.~$H_k(X)\in\QQ[X]$.

In this subsection, we would like to discuss the properties of these
class polynomials and their splitting modulo a prime $p$. This is used
twice in our algorithm: first with $p=2$, since we are going to start
from a curve in characteristic $2$; and second with a large odd $p$,
after the class polynomials have been computed, in order to build CM
curves over large finite fields.  By abuse of
notation, we also use $p$ to denote the prime ideal generated by the
rational prime $p$.
For simplicity, we restrict to $h_{K_0} = 1$, in which case we know
the number of isomorphism classes (see~Subsection~\ref{cmiso}).
Similar arguments will apply in the general case.

Let $A$ be an abelian surface with principal polarization $E$ of $CM$
type $(K,\Phi)$ with $j$-invariants $j_1$, $j_2$, $j_3$ and let $k_0$
be the field of moduli which is the unique subfield of $\CC$ with the
property:
An automorphism $\sigma$ of $\CC$ is the identity on $k_0$ if and only
if there exists an isomorphism $\lambda: (A,E)\rightarrow
(A^\sigma,E^\sigma)=(A,E)^\sigma$.
Obviously, we have $k_0:=\QQ(j_1,j_2,j_3)$.
Let $(K^*,\Psi)$ be the reflex type of $(K,\Phi)$.
We can characterize $k_0^* := k_0K^*$ in terms of class field theory.

\begin{theorem}[Main Theorem of Complex Multiplication,~\cite{ShimuraI}]
\label{haupttheorem}
Given a CM-type $(K,\Phi)$ with reflex type $(K^*,\psi)$.
Consider the ideal group $H_0$ of ideals $\aa$ in $K^*$ for which
there exists an element $\mu$ in $K$ such that
$$
\prod_j \psi_j(\aa) = (\mu) \mbox{ and }N(\aa) = \mu\conj{\mu}.
$$
The group $H_0$ contains the principal ideals, and the corresponding
unramified class field over $K^*$ is the field $k_0^*$.
\end{theorem}

For every CM type $(K,\Phi)$ we find $h_K$ isomorphism classes
of principally polarized abelian varieties (cf.~\cite{WengH})
and the polynomial
$$
G_k^\Phi(X)=\prod_{\sigma\in\Sigma^\Phi} (X-j_k^\sigma)
$$
(where $\Sigma^\Phi$ is the set of isomorphism classes of principally
polarized abelian surfaces with CM type $(K,\Phi)$) lies in
$K^*[X]$ by Theorem \ref{haupttheorem}. Since it is invariant under complex
conjugation, we even get $G_k^\Phi(X)\in K_0^*[X]$ where $K_0^*$ is the
real subfield of $K^*$. If $K$ is Galois, $H_k(X)=G_k^\Phi(X)$ and if $K$
is non-normal, $H_k(X)=G_k^{\Phi}(X)G_k^{\Psi}(X)$ where $\Phi$
and $\Psi$ are the two different CM types.
The polynomial $G_k^\Phi(X)$ does not need to be irreducible over $K_0^*$,
since $[k_0^*:K^*]$ can be smaller than $h_K$. More precisely, we have
$$
[k_0^*:K^*]=|I_{K^*}/H_0|=|I_{K^{**}}/H_K^{**}|\times U_0/U_1
$$
where $I_{K^*}$ is the ideal class group of $K^*$, $I_{K^{**}}$ is
the group of ideals in $I_K$ which are of the form
$\prod_j \psi_j(\cA^*)$ for some $\cA \in I_{K^*}$,
$H_{K^{**}}$ is the subgroup of principal ideals of $I_{K^{**}}$,
$U_0$ is the group of units in $K_0$ which are of the form
$N_{K/\QQ}(\cB)(\beta)^{-1}$ where
$\beta=N_{K/\QQ}(\cB)$ and $U_1$ the subgroup of units in
$K_0$ which are a norm of a unit in $K$ (see \cite{ShimuraI}, p.~112,
Note 3 and p. 114, Example 15.4 (3)).

If $h_K$ is odd, $K$ is non-normal and $N(\epsilon_0)=-1$, we can deduce
$[k_0^*:K^*]=h_K$ \cite{hecke, ShimuraI}. We expect
$G_k^\Phi(X)$ to be irreducible over $K_0^*[X]$ (in general this might
not be true, since $K^*(j_1,j_2,j_3)=k_0^*$ does not imply
$K^*(j_k)=k_0^*$ for a single $j_k$).

We now consider the abelian variety obtained by reducing the
invariants modulo a prime in $k_0^*$. Let
$p$ be a rational prime and $\gP\mid p$ be a prime ideal of
degree $f$ in $k_0^*$ such that $v_{\gP}(j_k)\ge 0$ for all
$k$. By reducing $j_k\bmod \gP$ we obtain a curve $C$ over
$\FF_{p^f}$. Let $J_C$ be its Jacobian.

Note that $v_{\gP}(j_k)<0$ only if the reduction of the
corresponding principally polarized abelian variety $A$ defined over a
number field $k\supseteq k_0^*$ modulo some prime $\qq$
lying above $\gP$ is superspecial, i.e.~$A\bmod\qq$
is isomorphic to the product of two supersingular elliptic
curves~\cite{goren}.

We will determine the $p$-rank of $J_C$. For this we use to following theorem.

\begin{theorem}
(\cite[Chapter~4~Theorem~1.1]{LangC},~\cite[Section 19]{ShimuraI})
\label{zeta1}
Let $\pi$ be the Frobenius endomorphism on the Jacobian $J_C$ obtained
by reducing $j_1$, $j_2$ and $j_3$ modulo $\gP$.
There exists an element $\pi_0$ in $K$ such $i(\pi_0)=\pi$ where $i$
denotes the embedding $\cO_K\rightarrow \End(A)$.
Moreover, with this $\pi_0$ one has $g(N_{k_0^*/K^*}(\gP)) =
\pi_0\cO_K$ where $g(\pp), \pp$ ideal in $K^*$, is the ideal $\aa$
in $\cO_K$ such that $\aa\cO_L =
\prod_{\psi_\alpha\in\Psi}\pp^{\psi_\alpha}\cO_L$ where $L$ is the
Galois closure of $K$.
\end{theorem}

The theorems so far allow us to determine the Frobenius endomorphism
for every prime ideal in $k_0^*$. We now consider the case that the
field $k_0=\QQ(j_1,j_2,j_3)$ does not contain $K^*$.

Let $G = \Gal(k_0^*|k_0)$. There exists an injective homorphism
$\pi:G\rightarrow \Aut(K)$ defined by
$i(\alpha)^\sigma=i(\alpha^{\pi(\sigma)})$.
The image $\pi(G)$ is a Galois group $\Gal(K|M)$ for some subfield $M$.
To determine the characteristic polynomial of the Frobenius endomorphism
over a smaller subfield we use the following theorem:

\begin{theorem}[\cite{LangC}, Chapter 4, Theorem 6.2]\label{zeta2}
Let $\pp_0$ be a prime in $k_0$ where $A$ has good reduction
and let $\pi_0$ be the corresponding Frobenius.
Let $T = \End(K/\CC)/\pi(G)$ be a set of representatives of
embeddings of $K$ into $\CC$ modulo $\pi(G)$.
Then for every $l\ne p$ the characteristic polynomial of the
Frobenius is given by
$$
\prod_{\tau\in T} \prod_{\gP\mid \pp_0}
(X^{f(\gP)}-\alpha(\gP^\tau))
$$
where $f(\gP)$ is the degree of the prime ideal $\gP$
in $k_0^*$ over $\pp_0$ where $\alpha(\gP)$
is up to a root of unity equal to $g(\pp_1)$ where
$\pp_1$ is the prime ideal in $K^*$ lying below $\gP$.
\end{theorem}

\begin{theorem}\label{lreihe}
Let $A$ be a principally polarized abelian surface with CM type
$(K,\Phi)$ with complex multiplication by $\cO_K$ with
invariants $j_1$, $j_2$ and $j_3$ in $k_0$. Let $(K^*,\Psi)$ be
the reflex CM type.
Consider the abelian variety $\redu{A}$ obtained by reducing
$j_k \bmod \pp_0$ for some ideal $\pp_0$ above $p$ in $k_0$.
Depending on the splitting of $p$ in $K$, we get
\begin{enumerate}
\item
if $p$ splits completely, the abelian variety
$\redu{A}$ is ordinary and has complex multiplication by
$\cO_K$;
\item
if $p$ is unramified, inert or splits only in $K_0/\QQ$ but not
any further, the abelian variety $\redu{A}$ is supersingular;
the same is true if $p$ ramifies completely, if $(p)=\pp^2$
and if $(p)=\pp_1^2\pp_2^2$ but $p$ does not
ramify in $K_0/\QQ$;
\item
if $p$ splits into three prime ideals, the abelian variety
will have $p$-rank 1; the same is true if
$(p)=\pp_1\pp_2\pp_3^2$;
\item
if $p$ is inert in $K_0/\QQ$ but splits in $K/K_0$, the abelian
variety will either be supersingular or ordinary with complex
multiplication by $\cO_K$ (depending on the CM type chosen);
the same happens if $(p)=\pp_1^2\pp_2^2$
where $p$ ramifies in the extension $K_0/\QQ$.
\end{enumerate}
\end{theorem}
\begin{proof}

\noindent
\begin{enumerate}
\item
Let $p \cO_K = \pp_1\overline{\pp_1}\pp_2\overline{\pp_2}$ with
all these prime ideals being distinct (since $p$ is unramified)
and let $\gP\mid \pp_0$ be the prime ideal in $k_0^*$. Then
$g(N_{k_0^*/K^*}(\gP))=(\pp_1\pp_2)^f=\mu\cO_K$
  is principal where $f$ is the degree of $N_(k_0^*/K^*)(\gP)\in K^*$ or
  equivalently the smallest integer such that
  $(\pp_1\pp_2)$ is principal.
We have $(\conj{\pp_1\pp_2)^f}$ is coprime to
  $(\pp_1\pp_2)^f$. Hence, the abelian variety is
  ordinary and by \cite{ShimuraI}, p.100, its endomorphism ring is
  equal to $\cO_K$.
\item If $p$ is inert, the abelian variety modulo $\pp_0$ is
  defined over $\FF_p$, $\FF_{p^2}$ or $\FF_{p^4}$ depending on
  whether $k_0\cap K^*=\QQ,K_0^*$ or $K^*$.
  Using Theorem~\ref{zeta1} and~\ref{zeta2}, we see that in this
  case the characteristic polynomial of the a power of the Frobenius is equal
  to $(X^2\pm p^4)^2$ and the abelian variety is supersingular.

  If $p$ splits in $K_0/\QQ$ but not any further, the abelian variety
  is defined over $\FF_p$ or $\FF_{p^2}$ and the characteristic
  polynomial of a power of the Frobenius is  $(X^2\pm p^2)^2$. Again, it is
  supersingular.

  If $(p)=\pp^4$, we have $g(N_{k_0^*/K^*}(\gP))$ is
  an even power of $\pp$. Considering the characteristic
  polynomial of the Frobenius we see that its $p$-rank must be equal
  to 0. The same argument works for $(p)=\pp^2$ and $(p)=\pp_1^2\pp_2^2$.

\item If $p$ splits in three prime ideals, then $p$ is inert in
  $K_0^*/\QQ$. This can only occur if $K$ is non-normal and in this case
  the field generated by the invariants will always contain $K_0^*$. Hence the
  field of definition will always contain $\FF_{p^2}$.

Let us consider the following diagrams of fields:
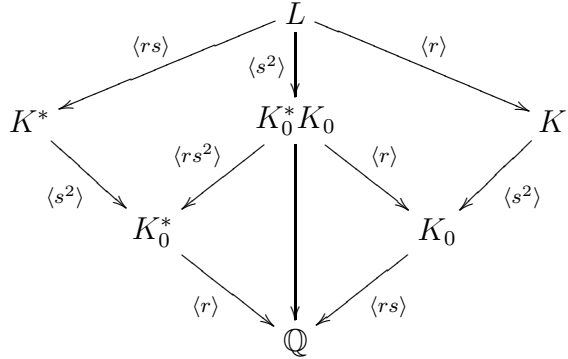
\begin{align*}
\xymatrix{
& & L \ar@{->}[ldl]_{\<rs\>}
      \ar@{->}[d]_{\<s^2\>}
      \ar@{->}[rdr]^{\<r\>} & & \\
K^\ast \ar@{->}[dr]_{\<s^2\>} & &
K_0^\ast K_0
      \ar@{->}[dl]_{\<rs^2\>}
      \ar@{->}[dd]
      \ar@{->}[dr]^{\<r\>} & & K \ar@{->}[dl]^{\<s^2\>}\\
 & K_0^\ast \ar@{->}[dr]_{\<r\>} & & K_0 \ar@{->}[dl]^{\<rs\>}&\\
 & & \QQ & &
}
\end{align*}
where $\<r,s;s^4=r^2=1,rs^3=sr\>$ is a representation of the
Galois group of the Galois closure $L$ of $K$.

Let $(p)=\pp_1\conj{\pp}_1\pp_2$ in $K$ and $\qq_1\qq_2
\conj{\qq}_1\conj{\qq}_2$ in $L$.
Then $(p) = \gR\conj{\gR}$, where $\gR =\qq_1 \qq_2$, is the prime ideal
decomposition of $p$ in $K^*$.

The automorphism $r$ leaves $\qq_1$, $\conj{\qq}_1$
invariant and interchanges $\qq_2$ and
$\conj{\qq}_2$. The automorphism $s^2$ maps $\qq_1$
to $\conj{\qq}_1$ and $\qq_2$ to $\conj{\qq}_2$.

The automorphism $r$ is a continuation of the real conjugation of
$K_0^*/\QQ$. We get
$\gR\gR^r=\qq_1\qq_2\qq_1^r\qq_2^r=\qq_1^2\qq_2\conj{\pp}_2$
Hence, $g(\gR)=\pp_1^2\pp_2$. The invariants are
defined over the field $\FF_{p^{2f}}$ where $f$ is the smallest
number such that $(\pp_1^2\pp_2)^f$ is principal.

In this case the Frobenius element is $(w)=(\pp_1^2\pp_2)^f$
and the Frobenius polynomial is the minimal polynomial of $w$.
By considering the Newton polygon of the characteristic polynomial of
the Frobenius we see that its $p$-rank is equal to $1$.
The case $(p)=\pp_1\pp_2\pp_3^2$ can be treated similarly.

\item
  Now consider the case where $p$ is inert in $K_0/\QQ$ but splits in
  $K/K_0$, i.e.~$p=\pp\conj{\pp}$. Here, $p$ splits in
  three prime ideals $(p)=\gR_1\conj{\gR}_1\gR_2$
  in $K^*$. We consider the same diagram as above.

We find $\gR_1\gR_1^r=\pp\cO_L$, hence
$g(\gR_1)=\pp\cO_K$, and
$\gR_2\gR_2^r=p\cO_L$, hence
$g(\gR_2)=p\cO_K$.

In the first case, the invariants are defined over the field
$\FF_{p^f}$ where $f$ is the smallest integer such that
$\pp^f$ is principal in $K$. The Frobenius element $w$ is then given
by $(w)=\pp^f$.

In the second case, the invariants are always defined over
$\FF_{p^2}$. We get $\pi_{0}=\pm p$. The abelian variety is
$\FF_{p^2}$-isogenous to the product of two supersingular elliptic
curves.
The case where $(p)=\pp_1^2\pp_2^2$ is similar.
\end{enumerate}
\end{proof}
Our algorithm will start with an ordinary hyperelliptic curve of genus
$2$ defined
over a finite field $\FF_{2^d}$. Hence, we are only interested in CM
field where $2$ splits completely, $2$ is inert in $K_0/\QQ$ but
splits in $K/K_0$ or $2$ ramifies in $K_0$ and is of the form
$\pp_1^2\pp_2^2$ in $K$.

We can then compute the extension degrees $d$ over which we expect to
find a hyperelliptic curve with complex multiplication by
$\cO_K$ as follows:
\begin{enumerate}
\item If $2$ splits completely, i.e.{}
  $(2) = \pp_1\conj{\pp}_1
         \pp_2\conj{\pp}_2$
  where $\conj{\pp}_i$ is the complex conjugate of
  $\pp_i$. Let $f_1$ be the smallest integer such that
  $\pp_1\pp_2$ is principal and let $f_2$ be
  the smallest integers such that
  $\pp_1\conj{\pp}_2$ is principal.
  Then we will find $h_K$ isomorphism classes of hyperelliptic
  curves defined over $\FF_{2^{f_1}}$ and $h_K$ isomorphism
  classes of hyperelliptic curves defined over $\FF_{2^{f_2}}.$
\item If $2$ is inert in $K_0/\QQ$ but splits into two prime ideals
  $\pp\conj{\pp}$ in $K/K_0$, we find $h_K$
  isomorphism classes of hyperelliptic curves with CM by
  $\cO_K$ over $\FF_{2^{f}}$ where $f$ is the smallest number
  such that $\pp^f$ is principal.
\item
If $2$ ramifies in $K_0/\QQ$ and is of the form $\pp_1^2\pp_2^2$
in $K/\QQ$, we find $h_K$ isomorphism classes of hyperelliptic curves
with CM by $\cO_K$ over $\FF_{2^{f}}$ where $f$ is the
smallest number such that $(\pp_1 \pp_2)^{2f}$ is principal.
\end{enumerate}

\section{The algorithm}
\label{main}

We now describe an algorithm for constructing hyperelliptic curves
over finite fields with complex multiplication by a given maximal order
$\cO_K$. We will restrict to specific CM-fields, e.g. there should
exist an ordinary hyperelliptic curve with complex multiplication
by $\cO_K$ over a field of characteristic 2 (see Subsection~\ref{split}).
The algorithm differs from the analytic approach mainly in the computation of
the class polynomials. Hence, we will first explain the construction
of the class polynomial.

\noindent
\textbf{Input:}
An ordinary hyperelliptic curve over $\FF_{2^d}$ with complex
multiplication by a maximal order $\cO_K$ in a CM field $K$.

\noindent
\textbf{Output:} Irreducible factors $\tilde{H}_k(X)$ of the
class polynomials $H_k(X)=\prod_{i=1}^s (X-j_k^{(i)})$ of degree
$n\le s$.
\begin{enumerate}
\item Compute the number of isomorphism classes $s$ of principally
  polarized abelian varieties over $\CC$ with complex multiplication by
  $\cO_K$ using Subsection \ref{cmiso}. This gives an upper
  bound for the degree of $\tilde{H_k}(X)$.
\item Lift the curve to a $2$-adic field.
\item Compute the Serre-Tate-Lubin lift using AGM.
\item Recover the absolute Igusa $j$-invariants $j_1$, $j_2$, $j_3$
  (see Section \ref{enum}) with high $p$-adic precision.
\item Apply $LLL$ to find the minimal polynomial $\tilde{H}_1(x)\in \QQ[x]$ of
  $j_1$ of degree $n\le s$.
\item  Apply $LLL$ to find the minimal polynomials $\tilde{H}_2(x)$ and $\tilde{H}_3(x)\in
  \QQ[x]$ of $j_2$ and $j_3$ of degree $n$.
\item Output $\tilde{H}_1$, $\tilde{H}_2$ and $\tilde{H}_3$.
\end{enumerate}
To get a curve over $\FF_p$ for $p$ a large prime with complex multiplication
by $\cO_K$, we now choose a prime $p$ such that there exists an
$w\in\cO_K$ with
$$
w\conj{w}=p
$$
and such that $f_w(1)$ is prime where $f_w$ is the minimal polynomial
of $w$.
We determine $(\redu{j}_1,\redu{j}_2,\redu{j}_3)\in\FF_p^3$
and compute the curve using Mestre's algorithm (cf.~\cite{WengA}).

\begin{remark} \label{remark_imp}

The determination of the number $s = s(K)$ of isomorphism classes is
useful for the application of the $LLL$ algorithm.
Note that $H_1(X)$ does not have to be irreducible
(cf.~Remark~\ref{remell}, (3)), but in many cases, we have
$H_1(X) = \tilde{H}_1(X)$ and using the algorithm above we can
recover the complete polynomial $H_1(X)$.
In general, we expect $n$ to be equal to $s$ or $s/2$.
This is true, for example, if $K$ is Galois, the real quadratic
subfield has class number one, the fundamental unit is negative,
and the class number $h_K$ is odd \cite{hecke}.
Hence, given $s$, we first try to take $n = s$ and $2n = s$
in step (5) of the algorithm.

Note that for most application (e.g. constructing curves over
large prime fields with given group order) it is sufficient to
compute an irreducible factor over $\QQ[x]$ of the class polynomial.
\end{remark}

\begin{remark}
\label{shortcut}
It is classical to consider that $\tilde{H}_1(X)$, $\tilde{H}_2(X)$ and
$\tilde{H}_3(X)$ are to be called the class polynomials. However, they do
not describe fully the CM points in the moduli space, since the relations
between the invariants are missing. We therefore modify the above
algorithm as follows.

Instead of computing $\tilde{H}_2(X)$ and $\tilde{H}_3(X)$ in (6),
we compute polynomials $G_2(X)$ and $G_3(X)$ of degree $n-1$ such that
$$
j_2\tilde{H_1}^\prime(j_1) = G_2(j_1)
\mbox{ and } j_3 \tilde{H_1}^\prime(j_1) = G_3(j_1)
$$
(see in Subsection \ref{sec:modlll} why this is better than classical
interpolation).\\
This approach is only possible if the coordinate $j_1$ is a separating
function for the CM points, or equivalently if $\tilde{H}_1(X)$ is a
squarefree polynomial of maximal degree. This is usually the case, and
then there is a major advantage for the second part of the
algorithm, the application of Mestre's algorithm~\cite{mestre}.

For Mestre's algorithm we reduce the polynomials modulo $p$ and
we try to find a suitable triple
$(\redu{j}_1,\redu{j}_2,\redu{j}_3)\in\FF_p^3$.
Given $\tilde{H}_1(x)$, $\tilde{H}_2(x)$ and $\tilde{H}_3(x)$
we have to loop through all possible triples $(x_1,x_2,x_3)$
where $\tilde{H}_k(x_i)=0$. In our situation we can compute
a root $\redu{j}_1$ of $\tilde{H}_1(X) \bmod p$ and then
determine $\redu{j}_2$ and $\redu{j}_3$ directly from
$G_2(x)$ and $G_3(x)$.
We need to factor only one polynomial modulo $p$ and the set
$(\redu{j}_1,\redu{j}_2,\redu{j}_3)\in\FF_p^3$
can be deduced directly. This is much more efficient.
Note that this trick can also be applied in the analytic approach.
\end{remark}

\section{Isomorphism classes in characteristic 2}
\label{enum}

In this section we discuss the choice of suitable invariants $j_1$,
$j_2$, $j_3$ for the algorithm described in Section \ref{main} and
describe how to choose suitable curves in characteristic 2 which we can use
as an input for the algorithm described in the previous section.

We have to be careful to choose the right invariants, since
we are in characteristic 2. In the literature and the computer
algebra package Magma we usually find three different sets of
invariants, known as Igusa-Clebsch invariants, Clebsch invariants
and Igusa invariants. These can easily be transformed into each
other (see \cite{mestre}).

We are only interested in the so-called Igusa invariants
$[J_2,J_4,J_6,J_8,J_{10}]$ since they also make sense in
characteristic 2.  Two curves are isomorphic over some
extension if and only if their Igusa invariants agree as
points in weighted projective space.  The subspace of
curves with ordinary Jacobian in $\bar{\FF}_2$ is defined
by the condition $J_2 \ne 0$~(see \cite{cnp}).
From the projective Igusa invariants we define absolute
invariants
$$
(j_1,\,j_2,\,j_3,\,j_4,\,j_5) =
(\frac{J_2^5}{J_{10}},
 \frac{J_2^3J_4}{J_{10}},
 \frac{J_2^2J_6}{J_{10}},
 \frac{J_2J_8}{J_{10}},
 \frac{J_4J_6}{J_{10}}).
$$
The absolute invariants are well-defined since $J_{10} = 0$
if and only if a curve is singular.  However, since $4J_8 =
J_4^2 - J_2J_6$, in characteristic~$2$ we obtain the relation
$j_1j_3 = j_2^2$ and also $j_2j_3 = j_1j_5$.
For an ordinary curve, the absolute invariant $j_1$ is nonzero
(since $J_2\ne0$) so we may eliminate the invariant $j_3$
in determining a parametrization of such curves, but use
the triple $(j_1,j_2,j_4)$ for the invariants of a lifted
curve.  In order to classify curves of nonordinary Jacobian,
it is necessary to define additional absolute invariants
(see~\cite{igusa}).

We now would like to enumerate all isomorphism classes of hyperelliptic
genus~2 curves over $\clos{k}$ which are defined over $k=\FF_{2^d}$
to find suitable CM fields as input to our algorithm.  We note that
over a finite field, a curve is defined over its field of moduli,
hence the field of definition of the point $(j_1,j_2,j_3,j_4,j_5)$ is
the field of definition for a curve.  (For a classification of curves
and their twists, we refer to a paper by Cardona, Nart and
Pupol{\'a}s~\cite{cnp}).

Following Igusa~\cite{igusa}, every ordinary curve of genus~2
in characteristic~2 has a normal form
$$
y^2 - y = ax + bx^{-1} + c(x-1)^{-1}, \quad abc \ne 0,
$$
isomorphic via $(x,y) \mapsto (x,y(x(x-1))^{-1})$ to the curve
$$
C : y^2 - x(x-1)y = x(x-1)(ax^3 + ax^2 + (b+c)x + b).
$$
We define $s_1(C) = a+b+c$, $s_2(C) = ab+bc+ac$, and
$s_3(C) = abc$.  The absolute Igusa invariants can be expressed
in terms of these invariants (cf.~\cite[p.\,623]{igusa});
in particular
$$
\begin{array}{l}
j_1^{-1} = J_2^{-5}J_{10} = s_3(C)^2,\\
j_2j_1^{-1} =
  J_2^{-2}J_4 = s_1(C)^2,\\
j_4j_1^{-1} =
  J_2^{-4}J_8 = s_2(C)^2 + s_1(C)^3 + s_1(C)^4.\\
\end{array}
$$
Thus the maps
$$
(s_1,s_2,s_3) \longmapsto \left(
  \frac{1}{s_3^2},\,
  \frac{s_1^2}{s_3^2},\,
  \frac{s_2^2 + s_1^6 + s_1^8}{s_3^2}\right),$$
and
$$
(j_1,j_2,j_4) \longmapsto \left(
  \frac{\sqrt{j_2}}{\sqrt{j_1}} ,\,
  \frac{\sqrt{j_4}}{\sqrt{j_1}} + \frac{j_2^2}{j_1^2} +
  \frac{j_2\sqrt{j_2}}{j_1\sqrt{j_1}},\,
  \frac{1}{\sqrt{j_1}}\right),
$$
define mutual inverses between triples $(s_1,s_2,s_3)$ with
$s_3 \ne 0$ and $(j_1,j_2,j_4)$ with $j_1 \ne 0$.

Conversely given a triple $(s_1,s_2,s_3) \in k^3$, with $s_3 \ne 0$,
there exists a curve in normal form
$$
C: y^2 - x(x-1)y = x(x-1)(ax^3 + ax^2 + (b+c)x + b).
$$
where $x^3+s_1x^2+s_2x+s_3 = (x-a)(x-b)(x-c)$, over an extension
of degree at most $3$.

\begin{remark}
Cardona, Nart and Pupol{\'a}s~\cite{cnp}) show for a finite field
of characteristic~$2$, that one can in fact find a representative
curve $C/k$ given any triple $(s_1,s_2,s_3)$ in $k^3$ with
$s_3\ne0$.
This implies that triples $(s_1,s_2,s_3)$ are in bijection with
$\overline{k}$-isomorphism classes of curves over $k$.  However, since
we require a curve in the split normal form above as input to our
algorithm of Section~\ref{main}, we omit the description of this
representative over $k$.
\end{remark}

\section{Higher dimensional generalization of AGM}
\label{agm2}

In order to apply the algorithm, we need a method to compute the
Serre-Tate-Lubin lift of a genus 2 curve over a $2$-adic field.
We use an algorithm due to Mestre \cite{Mestre00} which uses the explicit
formulae usually called ``Richelot isogeny''. Later on, Mestre
\cite{Mestre02} proposed another method, based on Borchardt's mean. The
latter has been implemented by Lercier and Lubicz \cite{LeLu03}.
Borchardt's mean involves simpler formulae and extends to higher genus.
Since we are interested in the genus 2 case, we stick to the first
``Richelot'' algorithm. This variant is not well described in the
literature, so we give now a few details about it.

\subsection{AGM lifting via Richelot's isogeny}

The basic idea of the genus 2 AGM lifting algorithm is to have explicit
formulae that describe fully a $(2,2)$-isogeny between jacobians of
curves. This can also be viewed as an explicit modular equation relating
the invariants of the curves. The following can be found in \cite{bm}:

\begin{theorem}
If $S$ and $T$ are monic polynomials of degree 2, define
$$ [S, T](x) = S'(x)T(x) - S(x)T'(x).$$
Let $C$ be a genus $2$ curve of equation $y^2=P(x)Q(x)R(x)$,
where $P$, $Q$, $R$ are monic of degree 2.  Let
$C'$ be the curve given by the equation
$$ \Delta y^2 = [Q, R](x)\ [R, P](x)\ [P, Q](x), $$
where $\Delta$ is the determinant of $P, Q, R$ in the basis $1, x,
x^2$.

Then $\Jac(C)$ and $\Jac(C')$ are $(2,2)$-isogenous abelian
varieties.  Moreover the kernel and the expression of the
isogeny can be made explicit.
\end{theorem}

This theorem is valid over any field of odd characteristic, including a
$2$-adic field. The next task is then to put the curve we have in a form
suitable to apply the theorem, and then to make the right choice for $P$,
$Q$ and $R$, so that the $(2,2)$-isogeny corresponds
to the second power
Frobenius isogeny, when we reduce everything modulo $2$.

A convenient form to work with is a Rosenhain form: we find $\lambda_0$,
$\lambda_1$ and $\lambda_\infty$ such that the curve of equation $y^2 =
x(x-1)(x-\lambda_0)(x-\lambda_1)(x-\lambda_\infty)$ is isomorphic to
$C$. By considering the reduction of the $2$-torsion divisors,
one can show that the $\lambda_i$ can be chosen such that
$\lambda_1\equiv 1 \bmod 4$, $\lambda_0\equiv 0 \bmod 4$ and
$\mathrm{val}(\lambda_\infty) = -2$.

Then the corresponding Rosenhain form for the curve $C'$, so that the
isogeny reduces to the second power Frobenius modulo 2, is given by
invariants $\lambda_i'$ satisfying
$$
\lambda_\infty' = \frac{(u_1-v_\infty)(u_\infty-v_0)}
{(u_1-v_0)(u_\infty-v_\infty)},\quad
\lambda_1' = \frac{(u_1-v_\infty)(w_1-v_0)}
{(u_1-v_0)(w_1-v_\infty)},\quad
\lambda_0' = \frac{(u_1-v_\infty)(w_0-v_0)}
{(u_1-v_0)(w_0-v_\infty)},
$$
where $u_1$ and $u_\infty$ are the solutions of the equation
$$ U^2 - 2\lambda_\infty U + \lambda_\infty
(1+\lambda_1)-\lambda_1 = 0,$$
$v_0$ and $v_\infty$ are the solutions of the equation
$$ V^2 - 2\lambda_\infty V + \lambda_0\lambda_\infty = 0,$$
and $w_0$ and $w_1$ are the solutions of the equation
$$ (\lambda_0-1-\lambda_1)W^2 + 2\lambda_1 W - \lambda_0\lambda_1 = 0.$$

In all these formulae, the subscript indicates the value of the variable
modulo $2$ (and an infinity subscript means that the valuation is
negative). Hence, the distinction between the roots of the equations of
degree 2 is easy.

As a consequence, we can derive a genus 2 AGM lifting procedure just like
in genus 1 as recalled in Section~\ref{genus1}. At each step, we have to
compute three square roots (for solving the three equations of degree 2) and a
few products, additions and inversions.  If the curve $C$ we started with is
ordinary, then the sequence converges (in the same sense as in
Section~\ref{genus1}) to the canonical lift of $C$. The theoretical
explanation for that is given in \cite{Carls04}.

To complete the algorithm, we still need to explain how to initialize the
AGM iteration. Since the formulae involve the 2-torsion points, we need
to have them defined in the base field that we consider. In other words,
when looking at the starting curve defined over the finite field $y^2 +
h(x) y = f(x)$, it is necessary that $h(x)$ splits completely. We
restrict to the case where $\deg h = 2$ and $\deg f = 5$. Also, since the
curve is supposed to be ordinary, the polynomial $h$ is squarefree. Let
us write $h(x) = x^2 + h_1x + h_0 = (x-\rho_0)(x-\rho_1)$. Then, by doing
the transformation to the Rosenhain form, and keeping everything formal,
we can derive the following values for the initialization of the AGM
iteration:
$$ \lambda_\infty = 4/h_1^2,\quad
\lambda_0 = 4\frac{f(\rho_0)h_1^2 + f'(\rho_0^2)}{h_1^6}, \quad
\lambda_1 = 1 + 4\frac{f(\rho_1)h_1^2 + f'(\rho_1^2)}{h_1^6}.$$

\subsection{Asymptotically fast lifting algorithm}

In the $p$-adic CM method, we might need to lift the curve to a very high
precision. The plain AGM method that we have just sketched has a
complexity which is at best quadratic in the precision. This quickly
becomes a problem. A first subquadratic algorithm was designed by Satoh,
Skjernaa and Taguchi \cite{SaSkTa03}, then an almost-linear lifting method
was designed by Kim et al. \cite{Kimetal02} in the case where the base
field admits a Gaussian normal basis, and finally Harley obtained an
almost-linear lifting method that works for any base field. A precise
description and comparison of these methods in the elliptic case can be
found in \cite{Vercauteren03}.

We have used the asymptotically fast variant of Harley, that we
now explain briefly.

Instead of going around the cycle of isogenous curves, getting closer and
closer to the canonical lift, we consider only two curves $C$ and $C'$
and their canonical lifts. Once lifted, the Rosenhain invariants of $C$
should annihilate the Frobenius-twisted modular equations
corresponding to the equations above: we should have
$$ \Phi(\Lambda, \Lambda^\sigma) = 0,$$
where $\Lambda$ is the vector $(\lambda_0, \lambda_1, \lambda_\infty)$ of
Rosenhain invariants, $\sigma$ is the  Frobenius substitution
in a $2$-adic field $\QQ_q$, and $\Phi$ is the function from $\QQ_q^6$
to $\QQ_q^3$
that corresponds to the Richelot equations above, where the intermediate
variables $u_i$, $v_i$ and $w_i$ have been eliminated.

Then an adaptation of the Newton lifting method can be used to compute a
solution $\Lambda$ to that equation, thus yielding the invariants of the
canonical lift. A key ingredient of that method is that we have to be
able to compute the action $\sigma$ quickly. To this effect, the $2$-adic
field is represented in a polynomial basis, with a generator that is a
root of unity (a Teichm\"uller lift of a generator of the underlying
finite field). Then the computation of the Frobenius image of an element
has a cost bounded by the cost of a few multiplications in the field.
We skip the details and refer to \cite{Vercauteren03} for a precise
description and analysis. Adapting the algorithm given there to the genus
2 case is essentially a multivariate rewriting of the algorithm for
elliptic curves. Not surprinsingly the jacobian matrix of $\Phi$ is
involved in place of just the two partial derivatives.

\section{The $p$-adic $LLL$ algorithm and Lagrange interpolation}
\label{lll}

As pointed out in Remark \ref{remell}, (3), the hyperelliptic curves
with complex multiplication by $\cO_K$ in characteristic 2
do not all have the same field of definition.
Moreover, given a class polynomial $H_K(X)\in\QQ[X]$, not all roots
in a field of characteristic $2$ need to lead to hyperelliptic curve
with complex multiplication by the field $K$
(see Subsection~\ref{split} and the discussion following
Theorem~\ref{lreihe}).
This happens for example for a non-normal $K$ whose real subfield has
class number one if $2$ is inert in the real subfield $K_0$ but splits
in $K/K_0$. In this case we find only $h_K$ hyperelliptic curves over
a field of characteristic 2 with complex multiplication by
$\cO_K$ although there exist $2h_K$ isomorphism classes
over $\CC$.

Hence, it is more convenient to compute only one root up to a high precision
and then apply the $LLL$ algorithm to recover the minimal polynomial. Note
that using this approach we will only find a irreducible factor of the class
polynomials and there are in general not irreducible.

\subsection{The $p$-adic $LLL$-algorithm}\label{lll1}
Given a lattice $\Lambda = \<b_1,\ldots,b_m\>$
the $LLL$ algorithm produces a short lattice basis.
This can be used to determine the minimal polynomial of an algebraic
element given by a floating point representation. Let $\det(\Lambda)$
be the determinant of $\Lambda$. Using Minkowski's inequality
we can approximate the shortest lattice vector by
$$
\sqrt{\frac{m}{2\pi e}}\det(L)^{1/m}.
$$
If $v\in \Lambda$ has length much smaller than this bound, it will
be the shortest vector with high probability.

Let $\ZZ_q$ be an extension of $\ZZ_2$ of degree d with $\ZZ_2$ basis
$1$, $w_1$,\ldots, $w_{d-1}$.  Let $\alpha \in \ZZ_q$ generating $\ZZ_q$,
and $\tilde{\alpha}$ be an approximation of $\alpha$ modulo a high
power of $2$, say $\alpha \equiv \tilde{\alpha} \bmod 2^N$.
We assume that we know the degree $n$ of its minimal polynomial
$f(x) \in \ZZ[x]$, i.e.
$$
f(x) = a_nx^n+\ldots+a_0
$$
where $a_i \in \ZZ$ are unknown.  In order to determine $a_i$, we
determine a basis of the left kernel in $\ZZ^{n+d+1}$ of the matrix
$$
\left(\begin{array}{c} A \\ 2^N I_d \end{array}\right)
$$
where $A$ is the $(n+1) \times d$ matrix
$$
\left(
\begin{array}{cccc}
1 & 0 & \ldots & 0\\
\alpha_{10} & \alpha_{11} & \ldots &\alpha_{1,(d-1)}\\
\vdots      &             &        &\vdots\\
\alpha_{n0} & \alpha_{n1} & \ldots &\alpha_{n,(d-1)}\\
\end{array}
\right)
$$
with $\alpha_{jk}$ defined by
$$
\alpha^j = \alpha_{j0}+\alpha_{j1}w_1+\ldots+\alpha_{j,(d-1)}w_{d-1}.
$$
This kernel is a lattice $\Lambda$, in which the coefficients of the minimal
polynomial of $\alpha$ are part of a short vector. Indeed, if
$a_0,\ldots,a_n$ are integers such that
$$
a_n\alpha^n+\ldots+a_0\equiv 0\bmod 2^N
$$
then $(a_0,\ldots,a_n,*,\ldots,*)$ will be a short vector in $\Lambda$         
that we expect to find in a $LLL$-reduced basis.                               

\subsection{Lagrange interpolation}
In Section \ref{main}, Remark \ref{shortcut}, we mention that
we do not compute $\tilde{H}_1(X)$, $\tilde{H}_2(X)$ and
$\tilde{H}_3(X)$ but $\tilde{H}_1(X)$ and two polynomials
$G_2(X)$, $G_3(X)$ with the property that
$$
j_2\cdot \tilde{H}_1^\prime(j_1)=G_2(j_1)\quad \text{and}\quad
j_3\cdot \tilde{H}_1^\prime(j_1)=G_3(j_1).
$$
Let us first consider the usual Lagrange interpolation, i.e.~suppose we
compute $F_k(X)\in\CC[X]$ with
$$
j_k=F_k(j_1) \mbox{ for $k=2,3$.}
$$
Let us assume that the conjugates $j_1^{(i)}$ for $i=1,\ldots,n$
are all distinct (see Remark \ref{distinct}).
Then $F_k(X)$ is given by
$$
\sum_{i=1}^{n} j_k^{(i)} \prod_{\ell\ne
i}\frac{X-j_1^{(\ell)}}{j_1^{(i)}-j_1^{(\ell)}}.
$$
Since $F_k(X)$ is easily seen to be Galois invariant, we have
$F_k(X)\in\QQ[X]$. Unfortunately, due to the factor
$$
\prod_{\ell\ne i}\frac{1}{j_1^{(i)}-j_1^{(\ell)}}
$$
the coefficients of $F_k(X)$ have usually a much larger height
than those of $\tilde{H}_k(X)$.  Hence, we prefer to compute
$G_k(X)$ with the property
\begin{equation}\label{gks}
j_k\tilde{H}_1^\prime(j_1)=G_k(j_1).
\end{equation}
A formula for $G_k$ is then given by
$$
G_k(X)=\sum_{i=1}^{n} j_k^{(i)} H^\prime(j_1^{(i)})
\prod_{\ell\ne i} \frac{X-j_1^{(\ell)}}{j_1^{(i)}-j_1^{(\ell)}}.
$$
Since
$$
H^\prime(j_1^{(i)})=\lc(\tilde{H}_1(X))\cdot \prod
(j_1^{(i)}-j_1^{(\ell)})
$$
where $\lc(\tilde{H}_1(X))$ denote the leading coefficient of
$\tilde{H}_1(X)$,
we expect $G_k(X)$ to have approximately the same height as $\tilde{H}_1(X)$.

\begin{remark}\label{distinct}
In order to be able to apply the Lagrange interpolation formula we need
the roots of the polynomial $\tilde{H}_1(X)$, to be distinct. In practice
we do not expect it to have any multiple roots. If this happens to be the
case, we solve the problem by choosing some linear combinations of $j_1$,
$j_2$, $j_3$ such that all roots are distinct.
\end{remark}

\subsection{Lagrange interpolation and $LLL$}
\label{sec:modlll}

We now modify the lattice given in Subsection~\ref{lll1}
to work for determining $G_2(X)$ and $G_3(X)$.
Let $d_k$ be the denominator of $G_k(X)$. Then equation
\eqref{gks} becomes
$$
d_kj_k^{(i)} \tilde{H}_1^\prime(j_1^{(i)})=\tilde{G_k}(j_1^{(i)})
$$
where $\tilde{G_k}(X)=d_kG_k(X)\in\ZZ[X]$.

We consider the lattice $\Lambda$ which is the kernel of the matrix
$$
\left(\begin{array}{c} A \\ 2^N I_d \end{array}\right)
$$
where the rows of $A$ contain the coefficients of $1$, $j_1$, $j_1^2$,
\ldots, $j_1^{n-1}$, $\tilde{H}_1^\prime(j_1)j_2$, expressed on the
$\ZZ_2$-basis.
If we have
$$
\tilde{G_k}(X)=a_{n-1}X^{n-1}+a_{n-2}X^{n-2}+\ldots+a_0,
$$
then the vector $(a_0,\ldots,a_{n-1},d_k,*,\ldots,*)$ will be a short
vector in $\Lambda$ that we expect to find in a $LLL$-reduced basis.

\subsection{Starting from several triples}\label{subsec:several}
It is often possible to compute $p$-adic approximations of several
triples  $(j_1, j_2,j_3)$ of invariants of curves having CM by
$\cO_K$. Furthermore, it can be the case that those triples form
an orbit under the action of a subgroup of the Galois group of the field
generated by the invariants.

Before showing how this can be used to speed up the computations, let us
give two examples of situations where we get such information.
\begin{itemize}
\item Once we have lifted one triple $(j_1,j_2,j_3)$ of elements in $\ZZ_q$ where
$q=2^d$, we can easily compute the $d$ conjugate triple by applying the
Frobenius automorphism of $\ZZ_q$.
\item It is possible that by enumerating
all isomorphism classes over the finite fields we have found several
nonconjugate curves having CM by $\cO_K$.  For instance, if $K$ is
non-normal, $h_{K_0}=1$, $N(\epsilon_0)=1$ and the class
number $h_K$ is odd, we expect to find at least $h_K$ isomorphism
classes over the finite field.
\end{itemize}

Let $(j_1^{(i)},j_2^{(i)},j_3^{(i)})_{1\leq i\leq k}$ be such a set of
conjugate triples, with $k$ divides $n$. Then the symmetric functions of
these triples are in an extension of degree $n/k$ of $\mathbb{Q}$. It is
then possible to build appropriate symmetric functions, so that applying
the $LLL$ algorithm to recognize algebraic numbers of degree $n/k$ will
allow to reconstruct the polynomials $H_1$, $G_2$, $G_3$. We expect this
approach to be faster than applying the $LLL$ algorithm to reconstruct
elements of degree $n$ directly, since the complexity of lattice
reduction depends badly on the dimension of the lattice (hence of the
degree of the elements to recognize).

On the other hand, having the possibility to recognize elements of
smaller degree implies more involved computations to deduce the
polynomials $H_1$, $G_2$, $G_3$. This is based essentially on resultant
computations. We give now more details about this approach.

We start by building the polynomial $M_1(X)$ whose roots are the
$j_1^{(i)}$:
$$
\begin{array}{rcl}
M_1(X) & = & (X-j_1^{(1)})(X-j_1^{(2)})\cdots (X-j_1^{(k)}) \\
       & = & X^k + m_{k-1}X^{k-1} + \cdots + m_1X + m_0.
\end{array}
$$
By the discussion above, the coefficients $m_i$ of $M(X)$ are algebraic
elements of degree $n/k$. We use the $LLL$ algorithm to compute the
minimal polynomial $P(X) \in \mathbb{Q}[X]$ of $m_0$. Let us call $K_P$
the number field $\mathbb{Q}[X] / (P(X))$, which is a degree $n/k$
subfield of the field $k_0$ of degree $n$ containing the CM invariants.
Then we recognize the other $m_i$ as elements of $K_P$, expressed in
terms of $m_0$. For that we use again the $LLL$ algorithm, but with the
modified lattice as in Section~\ref{sec:modlll}. Hence $M_1$ has been
rewritten as a bivariate polynomial
$$ X^k + m_{k-1}(Y)X^{k-1} + \cdots + m_1(Y)X + Y,$$
with rational coefficients, where $Y$ is a root of the the polynomial
$P(Y)$. The resultant in $Y$ of $M_1(X, Y)$ and $P(Y)$ is the polynomial
$H_1(X)$ we are looking for, perhaps up to a multiplicative factor.

We can perform the same kind of computation for $j_2$ and $j_3$, so as to
obtain $H_2(X)$ and $H_3(X)$. However, we would prefer to obtain $G_2(X)$ and
$G_3(X)$ that give more information. Let us explain how to get $G_2(X)$; the
polynomial $G_3(X)$ is computed in a similar manner.

Let $M_2(X)$ be the polynomial (with $p$-adic coefficients)
of degree at most $k-1$ such that $j_2^{(i)} = M_2(j_1^{(i)})$, for
$1\leq i\leq k$, that we can compute by a simple Lagrange interpolation.
Write $M_2(X) = n_{k-1}X^{k-1} + \cdots + n_1X + n_0$. As before, by
Galois invariance, the coefficients $n_i$ are algebraic elements of
degree $n/k$, and in fact are contained in $K_P$. We can recognize them
using the  $LLL$ algorithm with the modified lattice, and we get a
bivariate polynomial
$$ M_2(X, Y) = n_{k-1}(Y)X^{k-1} + \cdots + n_1(Y)X + n_0(Y),$$
defined over $\mathbb{Q}, $where $Y$ is again a root of $P(Y)$.
To convert back into a univariate representation, we need an explicit
expression for the embedding of the subfield $K_P$ into
$\mathbb{Q}[X]/(H_1(X))$. The computation of this embedding can be
handled by various algorithms. We suggest the following: the polynomial
$H_1(X)$ is obtained as the resultant of $M_1(X, Y)$ and $P(Y)$. If this
resultant is computed by the subresultant algorithm, on the way to the
solution we compute a polynomial of degree 1 in $Y$ that belongs to the
ideal generated by $M_1(X, Y)$ and $P(Y)$. Let us denote this polynomial
by $S(X,Y) = S_1(X) Y + S_0(X)$. Then as an element of
$\mathbb{Q}[X]/(H_1(X))$, a root of $P$ is given by $-S_0(X)/S_1(X)$,
thus yielding the required embedding.

Once $M_2$ has been recognized as an element of $\mathbb{Q}[X]/(H_1(X))$,
we just have to renormalize it with $H'_1(X)$, to obtained $G_2(X)$.

\begin{remark}
In the description of our method, we have overlooked two problems that
we encounter when actually implementing these algorithms:
\begin{itemize}
\item The elements are not algebraic integers, so we have to take care
of denominators everywhere. This is not a big difficulty but can induce
many programming mistakes.
\item If we implement line by line the method, there is a huge explosion
of the sizes of the coefficients in the middle of the algorithm. Once
$p$-adic elements are recognized as algebraic elements, we therefore have
to switch to modular computation: resultants, subresultants, and
computations in $\mathbb{Q}[X]/(H_1(X))$ must be handled by computing
modulo sufficiently enough primes, and we switch back to integers only
for the final reconstruction of $H_1$, $G_2$ and $G_3$, when we know that
the integers have a reasonnable size.
\end{itemize}
\end{remark}

\begin{remark}
As before, in this algorithm we made some genericity assumptions. Indeed,
it could well be that the coefficient $m_0$ that we used to defined the
subfield $K_P$ is in a fact in a subfield of degree less than $n/k$. In
that case, we just have to choose another element to define the field
$K_P$ we work with.
\end{remark}

\section{Determining the endomorphism ring}
\label{maxorder}

A  critical issue is the identification of a representative curve whose Jacobian has maximal
endomorphism ring.  It is necessary to have a mechanism to discard curves associated to the
nonmaximal orders. The following proposition gives a partial answer.

\begin{proposition}
\label{thm_maximal_order} Let $f$ be the minimal polynomial of the Frobenius endomorphism
 on the Jacobian $J_C$ of a genus $2$ curve $C$ defined over $\FF_q$ of
 characteristic $p$. Let $\pi$ be any root of this
polynomial and set $K = \QQ(\pi)$. Let the set
$$
\left\{
   \frac{g_1(\pi)}{p^{e_1} m_1}, \dots,
   \frac{g_t(\pi)}{p^{e_t} m_t} \right\}
$$
generate the maximal order $\cO_K$ over $\ZZ[\pi,\overline{\pi}]$ with $(m_i,p)=1$. Then
$g_i(\pi)/m_i$ is in $\End(J_C)$ if and only if $g_i(\pi)$ is the zero map on
$J_C[m_i](\overline{\FF_q})$.
\end{proposition}

\begin{remark}
If all the $e_i=0$, then we can really test the maximality, as
mentionned in \cite{EisLau04}. However,
and unlike the genus $1$ case, it is possible that
$\ZZ[\pi,\overline{\pi}]$ is not $p$-maximal in $\mathcal{O}_K$ and
then we cannot answer the problem.
\end{remark}
\noindent
Besides this algorithm, there are some other strategies which can be
applied :
\begin{enumerate}
\item Suppose we have given a curve $C$ of genus two with field of definition
$\FF_q=\FF_{2^d}$ and Frobenius polynomial $f_C(x)$. Let $K$ be the
quartic CM field generated by $f_C(x)$. Using the discussion following
Theorem \ref{lreihe}, we can compute the degree $f_1$ (resp. $f_2$) of
the field of definitions of the curves in characteristic 2 with complex
multiplication by $\cO_K$. If $d\ne f_1$ and $d\ne f_2$, the
endomorphism ring of $C$ cannot be maximal. Hence, we assume that
$d=f_i$ for some $i=1$, $2$.
\item Furthermore we can use the fact that the endomorphism ring of the
maximal order is in general as uncyclic as possible (a similar idea has been mentioned in
\cite{EisLau04}). By this we mean the following: Suppose we find two hyperelliptic curves $C_1$
and $C_2$ with the same characteristic polynomial i.e.~$f_{C_1}(x)=f_{C_2}(x).$ Then over every
field extension of $\FF_{2^d}$ the group of rational points on the Jacobian will have the same
order but not necessary the same group structure. Suppose we have a prime $\ell$ such that
$J_{C_1}$ has all $\ell$ torsion points rational ($\ell \ne p$) and not $J_{C_2}$, then the
conductor of the order of the endomorphism ring of $J_{C_2}$ will contain the prime $\ell$.
Indeed, $(\pi-1)/\ell \in J_{C_1}$ but is not in $J_{C_2}$.
\end{enumerate}
\section{Numerical examples}
\label{ex}

\subsection{Implementation}

We have implemented our algorithm using various computer algebra
packages. The first implementation has been written at a high level,
using the Magma system \cite{magma}. Then, to be able to deal with high
precisions, the asymptotically fast lifting algorithm using Richelot
isogeny has been implemented in C, based upon the {\tt Mploc} package
written by Emmanuel Thom\'e \cite{Mploc}. Finally, for the $LLL$
computations, we have interfaced our programs with Victor Shoup's NTL
library \cite{NTL}. Those three packages use the GMP library \cite{GMP}
for their time-critical integer operations.

After these optimizations, the cost of computing the canonical lift of a
curve is not so high, even if precision is huge. Therefore it appears
that the bottleneck of our method is the $LLL$ computation and the method
of section \ref{subsec:several} should be used for large examples.

\subsection{A non-Galois example with $n=2h$}
\label{nGalEx}

We start with the curve $C$ of equation $y^2 +h(x)y +f(x)=0$ over
$\mathbb{F}_8 = \mathbb{F}_2[t]/(t^3+t+1)$, with
$$ \begin{array}{rcl}
  f(x)  & = & x^5 + t^6x^3 + t^5x^2 + t^3x, \\
  h(x)  & = & x^2 + x.
  \end{array}
$$
The curve is ordinary and has CM by the maximal order of $K =
\mathbb{Q}(i\sqrt{23 + 4\sqrt{5}})$. The field $K$ is non-normal and its
class number is $3$; so we have $6$ isomorphism classes of principally
polarized abelian varieties.

We apply our algorithm and compute the canonical lift of $C$ to high
precision (in fact, a posteriori, we see that $1200$ bits are enough) and
get its invariants. From this we reconstruct the minimal polynomial $H_1$
and the corresponding $G_2$ and $G_3$. As expected, the degree of $H_1$
is $6$.
{\scriptsize
$$
\begin{array}{r@{\,}l}
H_1 & = 2^{18}5^{36}7^{24}\ T^6\\
&-\,11187730399273689774009740470140169672902905436515808105468750000\,T^5 \\
&+\,501512527690591679504420832767471421512684501403834547644662988263671875000\,T^4\\
&-\,10112409242787391786676284633730575047614543135572025667468221432704263857808262923\,T^3\\
&+\,118287000250588667564540744739406154398135978447792771928535541240797386992091828213521875\,T^2\\
&-\,2^{1}3^{50}5^{10}11^{1}13^{1}53^{1}701^{1}16319^{1}69938793494948953569198870004032131926868578084899317\,T\\
&+\,3^{60}5^{15}23^{5}409^{5}179364113^{5}\\
G_2  & =
2^{-3} \big(
2734249284974589542086559782016563911333032280921936035156250000\,T^5\\
&+\,57554607277149797568849387967258354564256002479144001401149377453125000000\,T^4\\
&+\,2402137816085408582966361480412923409977297040376760501014543382338189483861887923\,T^3\\
&-\,75691166837057576824962404339816428897154828109931810138346946500235981947587900092046875\,T^2\\
&+\,2^{1}3^{48}5^{10}35828519670812312117443096939126403484719666514876459782054400437\,T\\
&-\,3^{58}5^{15}11^{1}\,13^{2}23^{3}409^{3}23879^{1}179364113^{3}370974539856105277\big)\\
\multicolumn{2}{l}{G_3 = 2^{-4} \big(
200620022977265019387539624994933881234269211769104003906250000\,T^5}\\
&-\,23006467431764975697282545882188900514908468992554759536043135578125000000\,T^4\\
&+\,615017294619678068611319414718144161545088218260214211563850151291136646894987547\,T^3\\
&-\,14310698742415340178789612716269299249317950024503557714370659520249839645781463819312875\,T^2\\
&-\,2^{1}3^{46}5^{8}13^{1}61^{1}18373951326869^{1}25713288587261208212107985724468058651509734160907\,T\\
&+\,3^{55}5^{13}23^{2}409^{2}23561^{1}440131^{1}179364113^{2}451986402352017881724712641689\big)
\end{array}
$$
}
By looking at the Newton polygon of $H_1$ for the $2$-adic valuation, we
see that there are three roots that have valuation 0, and the others have
negative valuation. Hence only three of the curves have good reduction
modulo $2$. However, since $H_1$ is irreducible over $\QQ$, starting with
one curve (or from the 3 conjugate curves) yields the whole $H_1$.

This is consistent with Theorem~\ref{lreihe}. Indeed, $2$ is inert in
$K_0=\QQ(\sqrt{5})$ and splits in two prime ideals of degree $2$ in $K$.
Hence we are in subcase $(4)$. Furthermore, one can check that each of
the prime ideals above $2$ have order $3$ in the class group of $K$.

\subsection{A large example}

We start with the curve $C$ of equation $y^2 +h(x)y +f(x)=0$ over
$\FF_{32} = \mathbb{F}_2[t]/(t^5+t^2+1)$, with
$$ \begin{array}{rcl}
  f(x)  & = & x^5 + t^{20}x^3 + t^{17}x^2 + t^{19}x, \\
  h(x)  & = & x^2 + t^9x.
  \end{array}
$$
The curve is ordinary and has CM by the maximal order of $K =
\mathbb{Q}(i\sqrt{75 + 12\sqrt{17}})$. The field $K$ is non-normal and its
class number is $50$; so we have $s=100$ isomorphism classes of principally
polarized abelian varieties. The ideal $(2)$ splits completely in $K$, and
the primes above $2$ have order $5$ and $25$ in the class group.

However, when looking for a minimal
polynomial of the lifted value of $j_1$, the $LLL$ algorithm produced a
plausible answer of degree $50$. In fact, it seems that the class
polynomial of degree $s$ is not irreducible over the rationals, but
splits in two factors of degree $n=50$. Using our method, we can only
produce one of these factors $H_1(X)$, and the corresponding polynomials
$G_2(X)$ and $G_3(X)$.

For this large example, this would have been much faster to use the $5$
conjugate curves instead of only one. Indeed, with our implementation,
using only one curve (and therefore, doing lattice reduction to recognize
elements of degree $50$) requires about one day for the whole computation
on an Athlon64 processor, most of the time being spent in $LLL$.

For that case, we use a $p$-adic precision of $65000$ bits. The running
time to lift the curve and compute the invariants is $20$ seconds.

The leading coefficient of $H_1$ is
$3^{50}11^{156}17^{60}23^{72}41^{24}73^{12}83^{12}181^{48}691^{12}$.

\subsection{Checking the result}
Since we cannot give a bound on the coefficients of the class
polynomials, there is no way to  prove the result of the computation.
However there are some hints that indicate that the result is correct.
\begin{itemize}
\item The leading coefficient can be a large integer. However, we expect
this integer to be very smooth. In particular, it should be easy to
factor this number by trial division, even though the integer has several
hundreds of decimal digits. This could not occur for a random integer.
Therefore, if the answer of the $LLL$ algorithm has this property, then
we probably had enough precision.
\item When reducing the class polynomials modulo a suitable prime $p$,
one should be able to recover curves with the prescribed complex
multiplication. Hence, we can choose a prime $p$ small enough so that all
the computations are easy, and check that everything is consistent. For
instance, the large example of the previous section was checked with the
prime $p=47653$ which splits completely in $K$ into 4 prime ideals that
are principal. Then we check that $H_1$ splits completely over $\FF_p$,
and from its roots we deduce invariants and then equation for curves
(using Mestre's algorithm) that have indeed CM by $\cO_K$.
\end{itemize}

\section{complexity}
In this Section, we estimate the cost of our algorithm. The usual way of
computing class polynomials was described in \cite{WengA}. One
starts with a CM field, computes the period matrices
$(\Omega_i)_i$, recovers the $j$-invariants by computing theta
constants and computes the class polynomials by gathering all the
$j$-invariants. Weng's algorithm is dominated by the computation
of theta constants. This computation depends on the value of the
first minima of the period matrix, which makes the analysis of
this part difficult. However a naive evaluation of the theta
constants is quadratic in the precision. Our algorithm is linear
in the precision. Let us give some details. We can distinguish two
steps : the canonical lift of the curve and the LLL part.
Recursive programming based on the formulae of Richelot leads  to
a linear algorithm in the precision. More precisely the complexity
is $O((nk)^{1+\epsilon})$ where $n$ is the degree of the
extension, $k$ the final precision of the $p$-adic $j$-invariants
and $\epsilon$ represents the logarithmic factors in $n$ and $k$.
Then we use LLL to recover the class polynomials. Given $\langle b_1,
b_2, \ldots, b_m\rangle$ a basis of a lattice $\Lambda$ such that
for all $i$ in $\{ 1, \ldots, m \}$, $||b_i||^2 \leqslant B$, LLL
returns a LLL-reduced basis in a time $O(m^6\log^3(B))$. In our
case $\log(B)$ is the precision needed in order to make LLL work,
so this step is in $O(m^6k^3)$. The dimension of the
lattice $m$ is here the degree $h$ of our class polynomials
$\tilde{H}_k(X)$. Note that the floating-point version of LLL had
been improved by Nguy$\tilde{\mbox{\^e}}$n and Stehl{\'e} in \cite{NgSt05}. Their
version has a complexity of $O(m^5 (m+\log(B)) \log(B))$. When we
look at the LLL complexity, we can see that the dimension of the
lattice has a very bad influence on efficiency. To reduce the
dimension, one can proceed as suggested in  \ref{subsec:several}.
For instance, if one looks at the example \ref{nGalEx}, we can see
that $H_1(X)$ has  bad reduction modulo $2$. It gives a degree $3$
polynomial. Thus, if one seeks CM curves over $\FF_{2^3}$ with
maximal order in $K = \QQ(i\sqrt{23 + 4\sqrt{5}})$, one finds
three such curves. Hence, LLL has to deal with a lattice of
dimension only $2$. However such an enumeration is quite
expensive. It takes $2^{3n}$ operations to enumerate all the
curves and therefore one can afford it only over small extensions
of $\FF_2$. Note that, in practice, this idea is still valuable
because extensions of  $\FF_2$ of degree less than $10$ provide
already huge class number (for instance with $n = 7$, one can find
a quartic CM field whose class number is $6496$).

\section{Conclusion}

We have presented in this article a $2$-adic construction of CM
genus $2$ curves based on the AGM. This construction seems more
efficient than the existing complex method.
However, as for genus $1$, it does not allow to obtain all CM
fields. To tackle this problem, one should first find analogues of the
AGM method in characteristics greater than $2$.\\
Another possible generalization is to higher genus. Note that for
generic genus $3$ curves, any explicit  method is  known
to construct a curve over $\overline{\QQ}$ whose Jacobian has
complex multiplication. Such a construction can be done over
the $2$-adics with the AGM. However unlike the hyperelliptic case, one
does not know a complete set of invariants for non hyperelliptic genus
$3$ curves which, for the moment, prevent to make the link with number fields.

\bibliography{pcm2}

\newcommand{\etalchar}[1]{$^{#1}$}
\begin{thebibliography}{KPC{\etalchar{+}}02}

\bibitem[AM93]{AtkinB}
A.O.L. Atkin and F.~Morain.
\newblock Elliptic curves and primality proving.
\newblock {\em Math. Comp.}, 61:29--68, 1993.

\bibitem[Atk91]{AtkinA}
A.O.L. Atkin.
\newblock The number of points on an elliptic curve modulo a prime.
\newblock \textsl{{U}npublished manuscript}, 1991.

\bibitem[BC97]{magma}
W.~Bosma and J.~Cannon.
\newblock {\em Handbook of {M}agma functions}, 1997.
\newblock \\ {\tt http://www.maths.usyd.edu.au:8000/u/magma/}.

\bibitem[BLS02]{barreto}
P.~Barreto, B.~Lynn, and M.~Scott.
\newblock Constructing elliptic curves with prescribed embedding degrees.
\newblock {\em Security in Communication Networks -- SCN'2002, LNCS},
  2576:263--273, 2002.

\bibitem[BM88]{bm}
J.-B. Bost and J.-F. Mestre.
\newblock Moyenne arithm{\'e}tico-g{\'e}ometrique et p{\'e}riodes de courbes de
  genre 1 et 2.
\newblock {\em Gaz. Math. Soc. France}, 38:36--64, 1988.

\bibitem[BS04]{BS}
R.~Br\"oker and P.~Stevenhagen.
\newblock Elliptic curves with a given number of points.
\newblock {\em ANTS 2004, LNCS}, 3076:117--131, 2004.

\bibitem[BW03]{BW}
F.~Brezing and A.~Weng.
\newblock Elliptic curves suitable for pairing based cryptography.
\newblock {\em preprint}, 2003.

\bibitem[Car02]{Carls}
R.~Carls.
\newblock Mestre's method for point counting on elliptic curves in char $>$ 2.
\newblock {T}alk at Edchina, Sydney, 2002.

\bibitem[Car04]{Carls04}
R.~Carls.
\newblock {\em A generalized arithmetic geometric mean}.
\newblock PhD thesis, Rijksuniversiteit Gronigen, 2004.

\bibitem[CH02]{CH}
J.-M. Couveignes and T.~Hencoq.
\newblock Action of modular correspondences around cm points.
\newblock {\em ANTS V, LNCS}, 2369:234--243, 2002.

\bibitem[CNP04]{cnp}
G.~Cardona, E.~Nart, and J.~Pupol{\'a}s.
\newblock Curves of genus two over fields of even characteristic.
\newblock {\em to appear in Math. Zeitschrift}, 2004.

\bibitem[Coh96]{cohen}
H.~Cohen.
\newblock {\em A course in Computational Algebraic Number Theory}.
\newblock Springer, 1996.

\bibitem[DEM04]{dem}
R.~Dupont, A.~Enge, and F.~Morain.
\newblock Building curves with arbitrary small {MOV} degree over finite fields.
\newblock {\em to appear in Journal of Cryptology}, 2004.

\bibitem[dSG97]{goren}
E.~de~Shalit and E.~Z. Goren.
\newblock On special values of theta functions of genus two.
\newblock {\em Ann. Inst. Fourier (Grenoble)}, 47:775--799, 1997.

\bibitem[EL04]{EisLau04}
K.~Eisentraeger and K.~Lauter.
\newblock Computing {I}gusa class polynomials via the {C}hinese remainder
  theorem.
\newblock (http://www.arxiv.org/math.NT/0405305), 2004.

\bibitem[Gra02]{GMP}
T.~Granlund.
\newblock {\em The {GNU} {M}ultiple {P}recision arithmetic library -- 4.1}.
\newblock Swox AB, 2002.
\newblock distributed at \verb+http://swox.com/gmp/+.

\bibitem[Hec13]{hecke}
E.~Hecke.
\newblock {\"U}ber die {K}onstruktion relativ abelscher {Z}ahlk\"orper durch
  {M}odulfunktionen in zwei {V}ariablen.
\newblock {\em Math. Ann.}, 74:465--510, 1913.

\bibitem[Igu60]{igusa}
J.-I. Igusa.
\newblock The arithmetic variety of genus two.
\newblock {\em Ann. Math.}, 72:612--649, 1960.

\bibitem[Koh96]{kohel}
D.~Kohel.
\newblock {\em Endomorphisms of elliptic curves over finite fields}.
\newblock PhD thesis, University of California, Berkeley, 1996.

\bibitem[Koh03]{AGMX0N}
D.~R. Kohel.
\newblock The {AGM}-${X}_0({N})$ {H}eegner point lifting algorithm and elliptic
  curve point counting.
\newblock {\em ASIACRYPT 2003, LNCS}, 2894:124--136, 2003.

\bibitem[KPC{\etalchar{+}}02]{Kimetal02}
H.~Kim, J.~Park, J.~Cheon, J.~Park, J.~Kim, and S.~Hahn.
\newblock Fast elliptic curve point counting using {G}aussian normal basis.
\newblock In C.~Fieker and D.~R. Kohel, editors, {\em ANTS-V}, volume 2369 of
  {\em Lecture Notes in Comput. Sci.}, pages 292--307. Springer--Verlag, 2002.

\bibitem[KW04]{WengC}
K.~Koike and A.~Weng.
\newblock Constructing {CM}-{P}icard curves for cryptography.
\newblock {\em to appear in Math. Comp.}, 2004.

\bibitem[Lan83]{LangC}
S.~Lang.
\newblock {\em Complex Multiplication}.
\newblock Springer, 1983.

\bibitem[LL]{LeLu03}
R.~Lercier and D.~Lubicz.
\newblock A quasi quadratic time algorithm for hyperelliptic curve point
  counting.
\newblock Preprint.

\bibitem[LRD04]{LerRibD04}
R.~Lercier and E.~Riboulet-Deyris.
\newblock Elliptic curves with complex multiplication.
\newblock Number Theory List {\tt <NMBRTHRY@LISTSERV.NODAK.EDU>}, 2004.

\bibitem[LST64]{lst}
J.~Lubin, J.-P. Serre, and J.~Tate.
\newblock Elliptic curves and formal groups.
\newblock {\em Mimeographed notes, available under
  http://www.ma.utexas.edu/users/voloch/lst.-html}, 1964.

\bibitem[Mesa]{Mestre02}
J.-F. Mestre.
\newblock Algorithmes pour compter des points de courbes en petite
  caractéristique et petit genre.
\newblock Talk given in Rennes in March 2002, notes taken by D. Lubicz.

\bibitem[Mesb]{Mestre00}
J.-F. Mestre.
\newblock Utilisation de l'{AGM} pour le calcul de {$E(F_{2^n})$}.
\newblock Lettre adress{\'e}e {\`a} Gaudry et Harley, D{\'e}cembre 2000.

\bibitem[Mes91]{mestre}
J.-F. Mestre.
\newblock Construction des courbes de genre 2 {\`a} partir de leurs modules.
\newblock {\em Prog.Math., Birkh\"auser}, 94:313--334, 1991.

\bibitem[NS05]{NgSt05}
P.~Nguyen and D.~Stehl{\'e}.
\newblock Floating-point {L}{L}{L} revisited.
\newblock To appear in Eurocrypt'05 proceedings, 2005.

\bibitem[Shi98]{ShimuraI}
G.~Shimura.
\newblock {\em Abelian Varieties with complex multiplication and modular
  functions}.
\newblock Princeton University Press, revised edition, 1998.

\bibitem[Sho]{NTL}
V.~Shoup.
\newblock {\em {NTL}: A library for doing number theory}.
\newblock distributed at \verb+http://www.shoup.net/ntl/+.

\bibitem[Spa94]{Spall1}
A.-M. Spallek.
\newblock {\em Kurven vom Geschlecht 2 und ihre Anwendung in
  Public-Key-Kryptosystemen}.
\newblock PhD thesis, Institut f\"ur Experimentelle Mathematik, Universit\"at
  GH Essen, 1994.

\bibitem[SST03]{SaSkTa03}
T.~Satoh, B.~Skjernaa, and Y.~Taguchi.
\newblock Fast computation of canonical lifts of elliptic curves and its
  application to point counting.
\newblock {\em Finite Fields and Their Applications}, 9:89--101, 2003.

\bibitem[Tho]{Mploc}
E.~Thom{\'e}.
\newblock {\em Mploc}.
\newblock A library for local field.

\bibitem[Ver03]{Vercauteren03}
F.~Vercauteren.
\newblock {\em Computing zeta functions of curves over finite fields}.
\newblock PhD thesis, Katholieke Universiteit Leuven, 2003.

\bibitem[Wat69]{Waterhouse}
W.C Waterhouse.
\newblock Abelian varieties over finite fields.
\newblock {\em Ann. Sci. École Norm. Sup.}, 2(4):521--560, 1969.

\bibitem[Wen01]{WengB}
A.~Weng.
\newblock Hyperelliptic {CM}-curves of genus 3.
\newblock {\em Journal of the Ramanujan Mathematical Society 16}, 4:339--372,
  2001.

\bibitem[Wen03]{WengA}
A.~Weng.
\newblock Constructing hyperelliptic curves of genus 2 suitable for
  cryptography.
\newblock {\em Math. Comp.}, 72:435--458, 2003.

\bibitem[Wen04]{WengH}
A.~Weng.
\newblock Improvements and extensions of the {CM} method for genus two.
\newblock {\em High Primes and Misdemeanours, Fields Institute Communications,
  Series Volume}, 41:379--389, 2004.

\end{thebibliography}

\end{document}